\documentclass[11pt]{article}
\usepackage{authblk}

\usepackage{geometry}
\usepackage{amsmath,amssymb,mathrsfs}
\usepackage{color}
\usepackage{graphicx}
\usepackage{subfigure}
\usepackage{soul}
\usepackage{bbding}
\usepackage{lineno,hyperref}
\usepackage{multirow, booktabs, stmaryrd}
\modulolinenumbers[5]

\DeclareMathOperator*{\argmin}{arg\,min}

\DeclareMathOperator*{\Loss}{Loss}

\def\bp{\mathbf{p}}
\def\bx{\mathbf{x}}
\def\bX{\mathbf{X}}
\def\bs{\mathbf{s}}

\title{A Shallow Ritz Method for Elliptic Problems with Singular Sources}

\author[1]{Ming-Chih Lai}
\author[1]{Che-Chia Chang}
\author[1]{Wei-Syuan Lin}
\author[2,3]{Wei-Fan Hu}
\author[1,3]{Te-Sheng Lin}

\affil[1]{Department of Applied Mathematics, National Yang Ming Chiao Tung University, Hsinchu 30010, Taiwan}
\affil[2]{Department of Mathematics, National Central University, Taoyuan 32001, Taiwan}
\affil[3]{National Center for Theoretical Sciences, National Taiwan University, Taipei 10617, Taiwan}

\begin{document}

\maketitle

\begin{abstract}
In this paper, a shallow Ritz-type neural network for solving elliptic equations with delta function singular sources on an interface is developed. There are three novel features in the present work; namely, (i) the delta function singularity is naturally removed, (ii) level set function is introduced as a feature input, (iii) it is completely shallow, comprising only one hidden layer. We first introduce the energy functional of the problem and then transform the contribution of singular sources to a regular surface integral along the interface. In such a way, the delta function singularity can be naturally removed without introducing a discrete one that is commonly used in traditional regularization methods, such as the well-known immersed boundary method. The original problem is then reformulated as a minimization problem. We propose a shallow Ritz-type neural network with one hidden layer to approximate the global minimizer of the energy functional. As a result, the network is trained by minimizing the loss function that is a discrete version of the energy. In addition, we include the level set function of the interface as a feature input of the network and find that it significantly improves the training efficiency and accuracy. We perform a series of numerical tests to show the accuracy of the present method and its capability for problems in irregular domains and higher dimensions.
\end{abstract}

%


\section{Introduction}

Fluid-structure interaction problems have many applications in science and engineering, one example of which is blood flow~(fluid) simulation in heart valve leaflets~(embedded structure)~\cite{peskin1977}. The numerical simulation challenge for such problems is mainly that the shapes of embedded structures are often irregular and such structures change with time. To address these issues, Peskin~\cite{peskin1977} proposed the so-called Immersed Boundary~(IB) method, which is widely used because of its simplicity of implementation. This method does not require body-fitted discretization of the structure, which can save significant computational efforts.

The immersed boundary method is, in fact, both a mathematical formulation and numerical method for fluid-structure interaction problems. In IB formulation, we represent the fluid variables in an Eulerian manner and the embedded structure in Lagrangian one. The embedded structure is usually one-dimensional lower than the fluid dimensional space and regards as a singular force generator. Thus, the governing equations comprise Navier-Stokes (or Stokes) equations with singular forces as the Dirac delta function.
To solve the problem numerically, a projection-type of Navier-Stokes solver (for instance, see an overview in \cite{GMS06}) often involves solving elliptic equations with singular sources for the intermediate velocity. The IB numerical method then solves the equations by using finite difference discretization and a smooth version of the discrete delta function to regularize the singular sources. This method is easy to implement but leads to first-order accuracy. Another example is solving heat equations with singular sources. If time discretization is applied (for instance implicit Euler method), then at each time step, an elliptic equation with singular sources must be solved. So motivated by the above applications, we aim to solve elliptic problems with singular sources on an interface in this paper.

As mentioned above, the IB method is first-order accurate for elliptic interface problems due to the discrete regularization of the delta function sources. There are several other grid-based methods that achieve better accuracy in literature. For instance, the immersed interface method (IIM) proposed by LeVeque and Li \cite{LL94} incorporates the jump conditions via local coordinates into the finite difference scheme so that the local truncation error near the interface can be first-order resulting in the overall second-order accuracy in maximum norm. A simple implementation version of IIM that directly uses the jump conditions developed by the first author and his coworkers in \cite{LT08, HLY15} has the second-order accuracy in maximum norm as well.
A boundary condition capturing method (also named as ghost fluid method (GFM)) proposed by Liu et al.~\cite{LFK00} is able to solve the elliptic interface problems in a dimension-by-dimension manner, and captures the solution and its normal derivative jumps sharply while smoothing the tangential derivative. The method is first-order accurate in maximum norm in general. Recently, Egan and Gibou~\cite{EG20} have developed a novel idea to extend the original GFM by recovering the convergence of the gradient so that second-order accuracy can be achieved without modifying the resultant linear system. Along with this GFM approach, many other numerical methods for solving elliptic problems with interfaces or in irregular domains \cite{GM12, GLTG15, EG17, BG19, GHF19, BG20} have been successfully developed to improve the overall accuracy in maximum norm. One should mention that the linear systems resulting from those above methods are often symmetric positive definite, which can be solved efficiently by using iterative methods. Nevertheless, like IIM, special numerical treatments are always needed in the finite difference discretization of GFM near the interface or the irregular domain boundary. Thus, a mesh-free neural network method for solving the above problems provides an alternative to circumvent the difficulty arising from discretization near the interface or domain boundary for grid-based methods.

Solving partial differential equations (PDEs) with deep neural networks (DNNs) has drawn much attention in the scientific computing community recently. Part of the theoretical reason can be attributed to the various kinds of expressive power for function approximations using DNN such as those described in \cite{Cybenko1989, HSW89, LPWHW17, HS18}, just to name a few.
In terms of implementation, there are mainly two different approaches; namely, the physics-informed neural networks~\cite{RPK19}, and the deep Ritz method~\cite{EY18}. The major difference between the two approaches is how the loss is defined. One trains the physics-informed neural networks by minimizing the mean squared error loss of the equation residual, along with the initial and boundary condition errors. The deep Ritz method, however, begins with formulating the variational problem equivalent to the original PDE, so the natural loss function in this framework is simply the energy. Both approaches share the same major mesh-free advantage and therefore can practically solve problems in complex geometry~\cite{SY21} and in high-dimensional space~\cite{HJE18, SS18}.

Regarding deep learning approaches to solve PDEs with solutions that are piecewise smooth, Wang et al.~\cite{WZ20} proposed a deep Ritz-type approach to solve elliptic interface problems with high-contrast discontinuous coefficients, and Cai et al.~\cite{CCLL20} introduced the deep least squares method to solve elliptic interface problems where the solutions are continuous but the derivatives have jumps across the interface. In~\cite{GY21}, a deep unfitted Nitsche method is developed where two deep neural networks are formulated to represent two components of the solution. In \cite{BLSZ20}, the authors developed an immersed boundary neural network for solving elliptic equations with singular sources in which the delta function singularity still appears in the formulation and is discretized by a smooth discrete delta function proposed by Peskin \cite{peskin1977}.

Recently, the authors of the present paper proposed an efficient and accurate Discontinuity Capturing Shallow Neural Network (DCSNN)~\cite{HLL21} to solve elliptic interface problems where the solution is only piecewise-continuous. By augmenting one coordinate variable, the proposed shallow neural network is trained with PINNs-type loss.

In this paper, we propose a new shallow Ritz method for solving elliptic problems with singular sources. The novelties of the proposed network are three-fold. First, we remove the delta function singularity appearing in the original PDE by formulating the variational problem. Second, we include the level set function, which is commonly used as an interface indicator, as an additional feature input of the network that effectively improves the model's efficiency and accuracy. Third, we approximate the solution using a shallow neural network with only one hidden layer that significantly reduces the training cost in contrast to DNN.

The rest of the paper is organized as follows. In Section~2, we show how to transform a $d$-dimensional elliptic equation with singular sources into an energy functional minimization problem.  The shallow Ritz method to solve the problem is presented in Section~3, followed by a series of numerical accuracy tests and comparisons in Section~4. We give some concluding remarks and future work in Section~5.

\section{Elliptic equations with singular sources on the interface}\label{sec:derivation}

We consider a $d$-dimensional elliptic equation in a bounded domain $\Omega \subset \mathbb{R}^d$, divided by an embedded interface $\Gamma \subset \mathbb{R}^{d-1}$ into two regions; namely, inside ($\Omega^-$) and outside ($\Omega^+$) of the interface, so that $\Omega = \Omega^- \cup \Omega^+ \cup \Gamma$. The interface $\Gamma$ is represented by its parametric form $\bX(\bs)$  with surface parametrization $\bs \in \mathbb{R}^{d-1}$. The elliptic equation with singular sources on the interface is written as
\begin{eqnarray}
\Delta  u(\bx) - \alpha u(\bx) & =  & f(\bx) + \int_\Gamma c(\bs) \delta(\bx - \bX(\bs))\mbox{ d}\bs, \quad \mbox{in }\Omega, \label{eq:poisson} \\
u(\bx)& = & g(\bx),  \quad \mbox{on } \partial\Omega, \label{eq:bdc}
\end{eqnarray}
where $\alpha$ is a non-negative constant ($\alpha=0$ corresponds to the Poisson equation), $f$ is a given function, $c(\bs)$ is the source density defined only on the interface, and $\delta$ is the $d$-dimensional Dirac delta function.
Throughout this paper, we only focus on the case of the Dirichlet boundary condition, but other boundary conditions can be easily applied without changing the major ingredients of the proposed method.

As mentioned before, Eq.~(\ref{eq:poisson}) appears in many realistic applications and, in particular, the fluid velocity equations in the immersed boundary formulation for fluid-structure interaction problems~\cite{Pes02}. The fundamental difficulty in numerically solving Eq.~(\ref{eq:poisson}) is that the second term of the right-hand side has the integral over the $d-1$ dimensional interface, while the delta function is $d$-dimensional, this leaves the term has one-dimensional delta function singularity.
The IB method uses finite difference discretization and a grid-based discrete delta function to approximate the integral term in Eq.~(\ref{eq:poisson}). Because of the singularity of the delta function, this regularization technique is only first-order accurate, no matter what discrete delta functions are used. It was rigorously proved by Li~\cite{Li15} that the IB method for Eqs.~(\ref{eq:poisson})-(\ref{eq:bdc}) in 2D is indeed first-order accurate.

It is known that the solution to Eq.~(\ref{eq:poisson}) is piecewise smooth across the interface; more precisely, the solution is continuous over the domain $\Omega$ but has a discontinuity in its normal derivative on the interface $\Gamma$. In fact, Eq.~(\ref{eq:poisson}) is equivalent to the following elliptic interface problem with the same boundary condition (\ref{eq:bdc}) as
\begin{align}\label{Eq:elliptic}
\Delta  u(\bx) - \alpha u(\bx) =  f(\bx) \quad \mbox{in }\Omega \setminus \Gamma, \quad \llbracket u\rrbracket(\bs)=0, \quad \llbracket\partial_n u\rrbracket(\bs)=c(\bs) \quad \mbox{on }\Gamma,
\end{align}
where the notation $\llbracket \cdot \rrbracket$ represents the difference of the quantity (from the outside value to the inside value) across the interface. The above derivation can be found for example in \cite{LI06}.
One can immediately see that the delta function singularity no longer exists in Eq.~(\ref{Eq:elliptic}), but in terms of normal derivative jump instead. It is also worth mentioning that using equation~(\ref{Eq:elliptic}), one can easily construct different analytic solutions for the purpose of numerical tests later.

As mentioned in the Introduction, we have developed a neural network solver called DCSNN~\cite{HLL21} to approximate piecewise continuous functions and to solve more general elliptic interface problems than the one in Eq.~(\ref{Eq:elliptic}). The DCSNN augments an additional coordinate variable to label the pieces of a function so it inherently represents a discontinuous function.
However, as you can see from Eq.~(\ref{Eq:elliptic}), the solution is continuous across the interface. So in this paper, we propose a shallow Ritz-type neural network to solve Eqs.~(\ref{eq:poisson})-(\ref{eq:bdc}) by augmenting a level set function (continuous) as a feature input so that the solution is continuous and the energy is the natural loss function of the neural network. To proceed, we first reformulate Eq.~(\ref{eq:poisson}) into a variational problem as follows.



\subsection{Variational problem}

As shown above, the solution to the elliptic equation with singular sources on the interface, Eq.~(\ref{eq:poisson}), is continuous but has discontinuous derivatives across the interface.
Since the solution is not classically smooth, we are going to use the usual Sobolev spaces $H_0^1(\Omega)$ and $H^1(\Omega)$ to define the solution space. We assume the right-hand side function $f \in L^2(\Omega)$ and the source strength $c \in L^2(\Gamma)$. In order to take the boundary condition (\ref{eq:bdc}) into account, followed the work in \cite{LS03}, we assume that the boundary data $g$ has a smooth extension $\tilde{g}$ to the bounded domain $\Omega$ so that $u-\tilde{g} \in H_0^1(\Omega)$ and $\tilde{g}|_{\partial \Omega}=g$.
So the solution space can be defined as $H^1_g = \{u \in H^1(\Omega): u-\tilde{g} \in H_0^1(\Omega)\}$.
We now formulate equation (\ref{eq:poisson}) into its weak form by simply multiplying the test function $v \in H_0^1(\Omega)$ in Eq.~(\ref{eq:poisson}). Using the Green's first identity and the fact that test function $v$ vanishes at the boundary $\partial \Omega$, we obtain the following weak formulation: find $u \in H^1_g$ such that
\begin{equation}
-\int_\Omega \nabla u(\bx) \cdot \nabla v(\bx) \mbox{ d}\bx - \alpha\int_\Omega u(\bx) v(\bx) \mbox{ d}\bx = \int_{\Omega} f(\bx) v(\bx) \mbox{ d}\bx + \int_{\Gamma} c(\bs)v(\bX(\bs))\mbox{ d}\bs, \label{weak}
\end{equation}
for all $v \in H_0^1(\Omega)$. Note that the second term on the right-hand side comes from the definition of the Dirac delta function as
\begin{eqnarray}
& & \int_{\Omega} v(\bx) \int_\Gamma c(\bs) \delta(\bx - \bX(\bs))\mbox{ d}\bs \,  \mbox{ d}\bx
\nonumber\\
&=&
\int_\Gamma c(\bs) \int_{\Omega}  v(\bx)  \delta(\bx - \bX(\bs))\mbox{ d}\bx \,  \mbox{ d}\bs
=  \int_{\Gamma} c(\bs)v(\bX(\bs))\mbox{ d}\bs.
\end{eqnarray}
One can also derive the weak formulation by using the equivalent equation (\ref{Eq:elliptic}) with jump conditions in which the Green's first identity is applied separately in $\Omega^+$ and $\Omega^-$, and then summing up together to get Eq.~(\ref{weak}).
The existence and uniqueness of the above weak solution can be found in \cite{LS03}.

To adopt a Ritz-type neural network to solve the problem, we now rewrite the above weak formulation (\ref{weak}) to its equivalent minimization problem as follows. Find $u \in H^1_g$ such that $L[u]\leq L[v]$ for all $v \in H^1_g$, where the energy functional reads
\begin{equation}
\label{eq:energy}
L[v] = \frac{1}{2}\int_{\Omega} |\nabla v(\bx)|^2 \mbox{ d}\bx + \frac{\alpha}{2}\int_{\Omega} v(\bx)^2 \mbox{ d}\bx+ \int_{\Omega} v(\bx)f(\bx)  \mbox{ d}\bx + \int_\Gamma c(\bs)v(\bX(\bs)) \mbox{ d}\bs.
\end{equation}
As a result, just like the weak formulation, the delta function singularity disappears completely in the energy and its contribution becomes a regular integral over the interface $\Gamma$. Thus, we do not need to handle the singularity problem arising from the original equation~(\ref{eq:poisson}).

\subsection{Boundary condition enforcement}

We aim to solve the problem by using neural networks as the solution representation.
However, while neural networks are expressive in function approximation, and one hidden layer neural network can be a universal approximator (see a review in \cite{Pin99}), the boundary condition requirement in $H^1_g$ still seems to be infeasible.
Thus, we simply relax the boundary condition (\ref{eq:bdc}) by adding a penalty term to the energy  (\ref{eq:energy}) that reflects the penalty effect if the boundary condition is not satisfied exactly. The energy functional is therefore modified as
\begin{eqnarray}
\tilde{L}[v] &=& \frac{1}{2}\int_{\Omega} |\nabla v(\bx)|^2\mbox{ d}\bx + \frac{\alpha}{2}\int_{\Omega} v(\bx)^2\mbox{ d}\bx+ \int_{\Omega} v(\bx)f(\bx)\mbox{ d}\bx \nonumber\\
& &+ \int_\Gamma c(\bs)v(\bX(\bs))\mbox{ d}\bs + \beta\int_{\partial\Omega} (v(\bx)-g(\bx))^2\mbox{ d}\bx, \label{eq:energy2}
\end{eqnarray}
where $\beta$ is some positive penalty constant.

Now, it is interesting to see what the global minimizer will be for such a modified energy functional. To proceed,
let us decompose a function $v \in H^1(\Omega)$ as a sum of two functions by writing $v=u+h$ (the choice of $u$ and $h$ will be clear later).
We have its energy
\begin{eqnarray*}
\tilde{L}[v] &=& \frac{1}{2}\int_{\Omega} |\nabla u+\nabla h|^2\mbox{ d}\bx + \frac{\alpha}{2}\int_{\Omega} (u+h)^2\mbox{ d}\bx+ \int_{\Omega} (u+h)f\mbox{ d}\bx \\
& &+ \int_\Gamma c(u+h)\mbox{ d}\bs + \beta\int_{\partial\Omega} (u+h-g)^2\mbox{ d}\bx\\
&=& \tilde{L}[u] + \frac{1}{2}\int_{\Omega} |\nabla h|^2 \mbox{ d}\bx  + \frac{\alpha}{2}\int_{\Omega} h^2 \mbox{ d}\bx+ \beta\int_{\partial\Omega} h^2\mbox{ d}\bx \\
& & + \int_{\Omega} \left(\nabla u\cdot \nabla h + \alpha u h + hf\right) \mbox{ d}\bx + \int_\Gamma ch\mbox{ d}\bs+ \beta\int_{\partial\Omega} 2h(u-g)\mbox{ d}\bx,
\end{eqnarray*}
where the term $\int_{\Omega} \nabla u\cdot \nabla h \mbox{ d}\bx$ can be rewritten, using the Green's first identity under the smoothness assumption of $u$ in $\Omega^+$ and $\Omega^-$ (i.e. $u \in C^2(\Omega^\pm)$), as
\begin{eqnarray*}
\int_{\Omega} \nabla u\cdot \nabla h \mbox{ d}\bx  &=&
\int_{\Omega^-} \nabla u \cdot\nabla h \mbox{ d}\bx + \int_{\Omega^+} \nabla u\cdot \nabla h \mbox{ d}\bx \\
&=& -\int_{\Omega\setminus \Gamma} h \Delta u \mbox{ d}\bx
- \int_{\Gamma} h\llbracket\partial_n u\rrbracket \mbox{ d}\bs
+ \int_{\partial\Omega} h\partial_n u \mbox{ d}\bx.
\end{eqnarray*}
Therefore we obtain
\begin{eqnarray*}
\tilde{L}[v] &=&  \tilde{L}[u] + \frac{1}{2}\int_{\Omega} |\nabla h|^2 \mbox{ d}\bx  + \frac{\alpha}{2}\int_{\Omega} h^2 \mbox{ d}\bx+ \beta\int_{\partial\Omega} h^2\mbox{ d}\bx \\
& &+ \int_{\Omega\setminus \Gamma} h (-\Delta u +\alpha u + f) \mbox{ d}\bx
+ \int_\Gamma \left(c-\llbracket\partial_n u\rrbracket\right)h \mbox{ d}\bs
+ \int_{\partial\Omega} h\left(\partial_n u + 2\beta(u-g)\right) \mbox{ d}\bx.
\end{eqnarray*}
That is, if we choose $u$ as the solution of Eq.~(\ref{Eq:elliptic}) but with a modified Robin-type boundary condition
\begin{align}\label{Eq:elliptic2}
u(\bx) + \frac{1}{2\beta}\partial_n u(\bx)= g(\bx), \quad \bx\in\partial\Omega,
\end{align}
we should have
\begin{equation}
\tilde{L}[v] =  \tilde{L}[u]+ \frac{1}{2}\int_{\Omega} |\nabla (v-u)|^2 \mbox{ d}\bx  + \frac{\alpha}{2}\int_{\Omega} (v-u)^2 \mbox{ d}\bx + \beta\int_{\partial\Omega} (v-u)^2\mbox{ d}\bx.
\end{equation}
So it is clear that the solution of Eq.~(\ref{Eq:elliptic}) (or equivalently Eq.~(\ref{eq:poisson})) with boundary condition (\ref{Eq:elliptic2}) is the global minimizer of the energy functional $\tilde{L}$.
We thus formally write the solution to the variational problem as
\begin{equation}
u = \argmin_{v\in H^1(\Omega)} \tilde{L}[v], \label{mini1}
\end{equation}
where $H^1(\Omega)$ is the set of trial functions that does not require any constraint at the domain boundary. It is important to mention that as the penalty constant $\beta$ becomes larger, the above minimizer should approach to the solution of Eq.~(\ref{Eq:elliptic}) with exact boundary condition (\ref{eq:bdc}) (In fact, as $\beta \rightarrow \infty$, Eq.~(\ref{Eq:elliptic2}) tends to be Eq.~(\ref{eq:bdc})).

\subsection{Level set function augmentation} \label{smooth}

As shown above, the minimizer of the energy functional (or the solution of Eq.~(\ref{Eq:elliptic})) is a continuous function, but has a jump discontinuity at its normal derivative on the interface.  Therefore, the trial functions must carry the same feature. We note that the construction of such a set of functions using neural net approximation requires additional efforts.
Here, inspired by DCSNN~\cite{HLL21}, where an augmented variable is introduced to categorize precisely the spatial coordinates into each sub-domain, we also introduce an augmented variable and require the function to be continuous throughout the whole domain.
More precisely, consider a level set function $\phi(\bx)$ such that the zero level set gives the position of the interface $\Gamma$, i.e., $\Gamma = \{\bx \in\mathbb{R}^d\,|\, \phi(\bx)=0\}$.
We define a function $U(\bx, z): \mathbb{R}^{d+1}\to  \mathbb{R}$ that satisfies
\begin{equation}
u(\bx) = U(\bx, \phi(\bx)), \quad \bx\in\Omega.
\end{equation}
Here the level set function $\phi(\bx)$ is considered as a feature input in the augmented variable $z$. We require both the level set function $\phi$ and the extension function $U$ being continuous, so that $u(\bx)$ is a continuous function. Although it looks like the derivative discontinuity on the interface is not considered, the augmented variable $\phi(\bx)$ somehow gives additional information to the function $U$ related to the interface. As we will see later in the numerical experiments, indeed the introduction of the augmented variable effectively improves the capability of neural networks in function approximation.

It is worth mentioning that we assume the solution to the problem takes the form: $u(\bx) = U(\bx, \phi(\bx))$. Therefore, if one instead uses an indicator function for the interface as the additional feature input, e.g., $\phi(\bx) = 1$ if $\bx\in\Omega^+$ and $\phi(\bx) = -1$ if $\bx\in\Omega^-$, the resulting function will be discontinuous, which essentially violates the assumption of a continuous solution that we consider here.

\subsection{Summary}

To summarize, we solve elliptic equations with singular sources on the interface by looking for a $d+1$ dimensional continuous function of the form $U(\bx, z)$. The target function is found by minimizing the energy functional $\tilde{L}[u]$, where $u(\bx)= U(\bx, \phi(\bx))$.
In the following we shall develop a neural network architecture to represent $U$ and a loss function to be used for model training.

\section{A shallow Ritz method}

We propose a shallow neural network to approximate the $(d+1)$-dimensional continuous function of the form $U(\bx, \phi(\bx))$ and serve as an ansatz for solving the minimization problem Eq.~(\ref{mini1}). Based on the universal approximation theory~\cite{Cybenko1989}, we hereby design a \emph{shallow}, \emph{feedforward}, \emph{fully-connected} neural network architecture, in which only one hidden layer is employed.
Let $N$ be the number of neurons used in the hidden layer, the approximation function (or output layer) under this network structure is explicitly expressed by
\begin{align}\label{Eq:SINN}
u(\bx)=U(\bx, \phi(\bx)) = W^{[2]}\sigma(W^{[1]}(\bx, \phi(\bx))^T+b^{[1]}) + b^{[2]},
\end{align}
where $W^{[1]}\in\mathbb{R}^{N\times(d+1)}$ and $W^{[2]}\in\mathbb{R}^{1 \times N}$ are the weight matrices, $b^{[1]}\in\mathbb{R}^{N}$ and $b^{[2]}\in\mathbb{R}$ are the bias vectors, and $\sigma$ is the activation function. One can easily see that the function $U$ is a linear combination of the activation functions and hence shares exactly the same smoothness as $\sigma$.
By collecting all the training parameters (including all the weights and biases) in a vector $\bp$, the total number of parameters in the network (i.e., dimension of $\bp$) is counted by $N_p = (d+3)N + 1$.

As the solution we are looking for is the global minimizer of the energy functional, it is therefore natural to consider the loss function in the training procedure using that energy. In this stage, we aim to learn the training parameters $\bp$ via minimizing the modified energy (\ref{eq:energy2}).
The loss function in the training procedure is thus defined in the framework of deep Ritz method~\cite{EY18} by replacing the integrals in Eq.~(\ref{eq:energy2}) using discrete quadrature rules.
That is, given training points in $\Omega$, along the embedding interface $\Gamma$, and on the domain boundary $\partial\Omega$, denoted by $\{\bx^i\}_{i=1}^M$, $\{\bx^j_{\Gamma}\}_{j=1}^{M_\Gamma}$ and $\{\bx^k_{\partial\Omega}\}_{k=1}^{M_b}$, respectively, we define the loss function as
\begin{align}\label{Eq:loss}
\begin{split}
\Loss(\bp) &= \frac{\mbox{Vol}(\Omega)}{M}\sum_{i=1}^{M} \left(\frac{1}{2} |\nabla u(\bx^{i})|^2 + \frac{\alpha}{2}u(\bx^{i})^2 +u(\bx^{i}) f(\bx^{i})\right) \\
&+ \frac{\mbox{Vol}(\Gamma)}{M_{\Gamma}}\sum_{j=1}^{M_{\Gamma}}c(\bx^j_\Gamma) u(\bx^j_\Gamma)
+\beta\,\frac{\mbox{Vol}(\partial\Omega)}{M_b}\sum_{k=1}^{M_b}\left(u(\bx^k_{\partial \Omega}) - g(\bx^k_{\partial \Omega})\right)^2,
\end{split}
\end{align}
where $\mbox{Vol}(\Omega)$ is the volume of $\Omega$ in $\mathbb{R}^d$, while $\mbox{Vol}(\Gamma)$ and $\mbox{Vol}(\partial\Omega)$ are volumes of $\Gamma$ and $\partial\Omega$, respectively,  in $\mathbb{R}^{d-1}$.
Here, for brevity, we suppress the notation of the parameters $\bp$ in the solution $u$ as one can see their dependence through the equation (\ref{Eq:SINN}).  Also note that, the evaluation of $u$ at the domain training point $\bx^{i}$ is realized through the relation $u(\bx^{i}) = U(\bx^{i}, \phi(\bx^{i}))$ so the  feature input of the network is $(\bx^{i}, \phi(\bx^{i}))$. Similarly, the evaluations of $u$ at the training points $\bx^j_\Gamma$ in $\Gamma$ and $\bx^k_{\partial \Omega}$ in $\partial \Omega$ are defined in the same manner.

\subsection{Selection of training points}

As one can see, the quadrature rule used to derive the loss function (\ref{Eq:loss}) is of the Monte Carlo type, where the integrals in modified energy Eq.~(\ref{eq:energy2}) are approximated by the mean of the samples. However, since the energy is presented as integrals, one may expect a better performance by evaluating these integrals using more accurate numerical quadrature rules. Thus, if the problem under consideration is defined in regular domains such as a rectangle or a circle, the first intuitive way to select the training points might be the quadrature nodes based on, e.g., midpoint or Gaussian quadrature rules, just to name a few. (Note, if Gaussian nodes are chosen as the training points, one should change the weights in approximating the energy functional.) But, we found that the training of neural networks heavily depends on ``accurate" evaluations of these integrals. The question is not which quadrature rule is chosen, but how to ensure accurate evaluation of the integral value.

In traditional scientific computing approach, a convergence test must be performed if we want to evaluate an integral to the desired accuracy but without a prior information about the number of quadrature nodes required. Therefore, if we fix the number of nodes to compute the integral, the result might not be reliable. In neural network approach, the parameters change at each iteration of the training process meaning that the function presented by the network differs from iteration to iteration. For the Ritz method, if the training points are fixed throughout the entire optimization process to minimize the energy functional, we can not guarantee the numerical integration accuracy.

To ensure a reliable estimation of the energy, ideally, we should check the convergence of the numerical quadrature at each training step which is not practical in reality. In this paper, we adopt the strategy proposed in \cite{EY18}. That is, we use Monte Carlo integration to evaluate the loss function Eq.~(\ref{Eq:loss}) but re-select training points at each iteration of the optimization process. The spirit of this idea is like the mini-batch gradient descent method. Imagine that our training dataset contains all the points in the domain $\Omega$, and at each iteration step, we train the model based on only one mini-batch of the entire dataset. The training process is ended if certain stopping criteria are satisfied. We find that such a procedure, even though it seems to be less accurate in evaluating the energy at first glance, is more stable in the sense that the loss will remain close to the true energy.

\subsection{Selection of Optimizer}

Even though the energy admits a global minimizer, as shown in the previous section, there may be local minimizers in the space of network parameters. On the other hand, since we re-select the training points in each iteration, the loss landscape will be different for each iteration. So it is not suitable to use traditional descent method such as the gradient descend method for minimization. In this paper, we exploit the Adaptive Moment Estimation (Adam) optimizer~\cite{KB17}, that can deal with sparse gradients and non-stationary objectives, which is well-suited to the current aim. It has been empirically found that Adam has never underperformed SGD in general~\cite{CSNLMD19}.

\section{Numerical results}

In this section, we perform several numerical tests in two, three, and six dimensions for elliptic problems with singular sources using the developed shallow Ritz method. In the following examples, we choose sigmoid as the activation function. After the training process is complete, we check the accuracy by computing the error between the neural network solution and the exact solution of the problem. To do this, we randomly choose $N_{test}$ testing points lying in $\Omega$ to compute the relative $L_\infty$ and $L_2$ errors as
\[
\frac{\|u_\mathcal{S}-u\|_\infty}{\|u\|_\infty}, \quad \frac{\|u_\mathcal{S}-u\|_2}{\|u\|_2},
\]
respectively, where
\begin{align*}
\|u\|_\infty = \max_{1\leq i \leq N_{test}}|u(\bx^i)|, \quad
\|u\|_2 = \sqrt{\frac{1}{N_{test}}\sum_{i=1}^{N_{test}}(u(\bx^i))^2}.
\end{align*}
Here $u$ is the exact solution of Eqs.~(\ref{eq:poisson})-(\ref{eq:bdc}) while $u_\mathcal{S}$ is the solution obtained from  the present model. In general, we set $N_{test}=100M$, where $M$ is the number of training points in $\Omega$. We always terminate the training procedure after $50000$ iterations with a fixed learning rate $5\times10^{-3}$.

Throughout all numerical experiments conducted here, the training time of the present network typically takes a few minutes on a desktop with 8-core Intel i7-10700k CPU and a Nvidia 2080Ti graphic card. Particularly, we found that the training time usually scales linearly with the number of training points but only weakly with the number of neurons. Regarding the choice of the number of neurons, we usually increase the number of neurons until the loss value (also the approximation of the energy) tends to converge. As you can see in our numerical results, just a moderate number of neurons (up to $40$) is generally sufficient to give predictive accuracy less than $1\%$ relative error in $L_2$ norm.

In the following numerical experiments, we will just list the analytic form of the exact solution in $\Omega$ and the level set function to represent the interface $\Gamma$. The resultant function $f(\bx)$ and the source density $c(\bs)$ in Eq.~(\ref{eq:poisson}) can be easily computed via the equation itself and the derivative normal jump condition $\llbracket\partial_n u\rrbracket=c(\bs)$, respectively.

\paragraph{\textbf{Example 1}}

In the first example, we choose a square domain $\Omega=[-1, 1]\times [-1, 1]$ with a circular interface that can be labeled by the zero-level set of the function $\phi(x,y) = x^2+y^2-0.5^2$. The exact solution is chosen as
\begin{align}
    u(x,y) =
    \left\{
    \begin{array}{ll}
        - \ln\left(x^2+y^2\right) & (x,y) \in\Omega^+,\\
        -\ln(0.5^2) & (x,y) \in \Omega^-.\\
    \end{array}\right.
\end{align}
We choose $\alpha=0$ in this example.

\subparagraph{\textbf{Fixed training points}}

At first, we would like to show the over-fitting problem of using training points that are fixed throughout the entire optimization process. Instead of using the Monte Carlo method to approximate the energy, we use the midpoint rule with fixed uniform grid points on domain $\Omega$, interface $\Gamma$, and domain boundary $\partial\Omega$, which is supposed to be more accurate in approximating the integrals. Notice that, the midpoint quadrature rule leads to exactly the same formula for the loss function as in Eq.~(\ref{Eq:loss}). More precisely, in the selection of training points, we use $40\times 40$ uniform quadrature nodes in the domain $\Omega$ as
\begin{equation}
(x_i, y_i) = \left(-1 + \left(i-1/2\right) \Delta x, -1 + \left(j-1/2\right) \Delta y\right), \quad i=1,\cdots, 40, \quad j=1, \cdots, 40,
\end{equation}
where $ \Delta x=\Delta y=\frac{1}{20}$ is the grid size. Both the interface and the domain boundary have $160$ grids points, which are also uniformly distributed. So overall we have $M=1600$, $M_\Gamma=M_b=160$ which gives $1920$ training points in total. The shallow neural network consists of one hidden layer with $30$ neurons, for a total number of $N_p = 151$ parameters to be trained. The penalty constant is set by $\beta=200$. We note that the number of training points is much larger than the number of parameters used in this case.

\begin{figure}[h]
\centering
\includegraphics[scale=0.7]{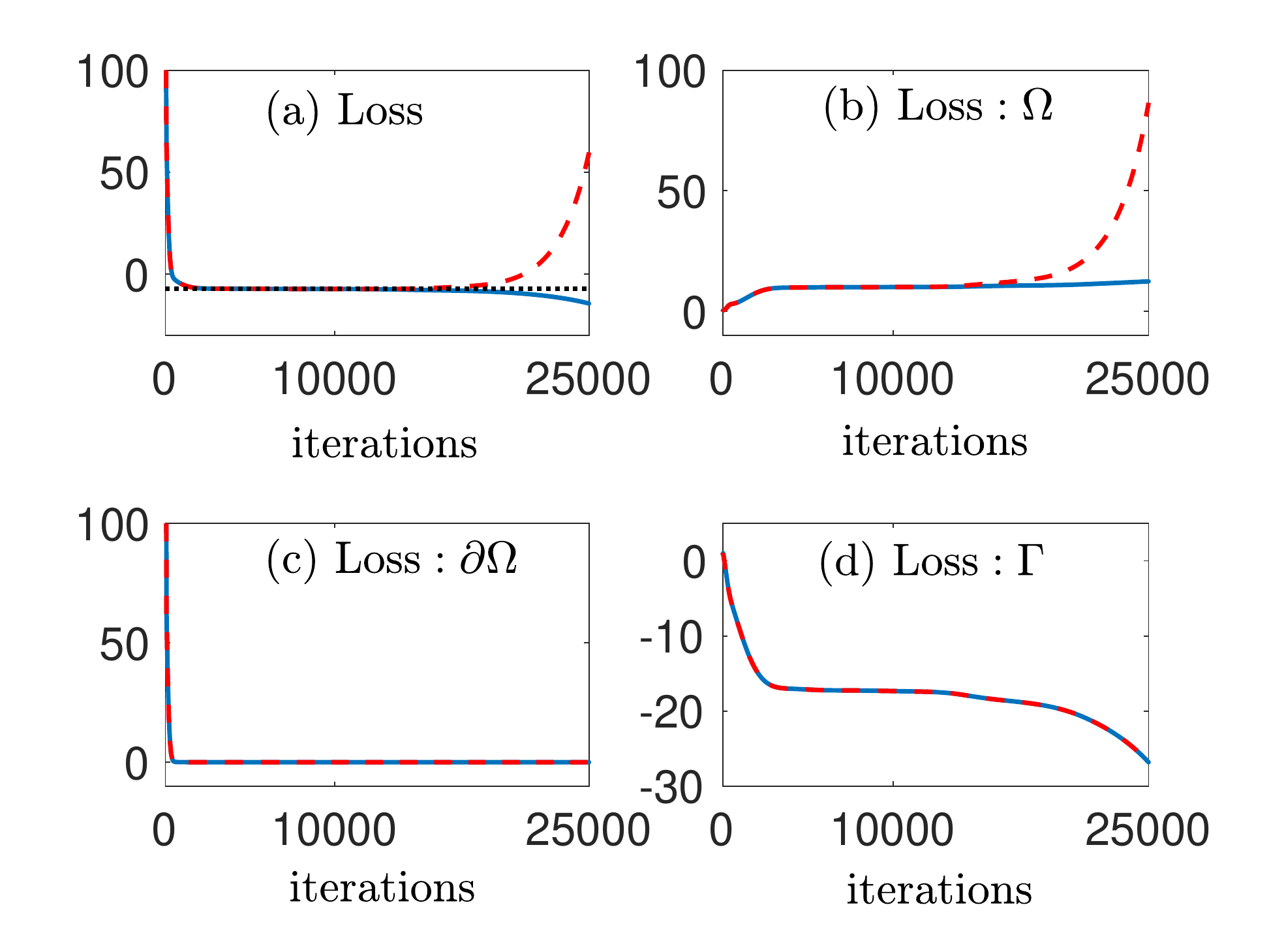}
\caption{Evaluation of the loss functions. The solid (blue) line shows the loss obtained during training process while the dashed (red) line shows a re-evaluation of the loss using testing points. The dotted (black) line in (a) indicates the theoretical energy value corresponding to the exact solution.}
\label{Fig:fixed_point_loss}
\end{figure}

The training loss obtained during the optimization process is shown in Fig.~\ref{Fig:fixed_point_loss}(a) as a solid (blue) line. One can see that the loss function indeed decreases throughout the entire training process. Meanwhile, the losses corresponding to the domain ($\Omega$), domain boundary ($\partial\Omega$) and the interface ($\Gamma$) are also shown in panels (b), (c) and (d) as solid (blue) lines, respectively.

In addition, in Fig.~\ref{Fig:fixed_point_loss}(a), we plot a dotted (black) line indicating the energy value that corresponds to the exact solution; that is, the theoretical minimum energy. As the loss function is a discrete representation of the energy, one might expect to have the same global energy minimizer so the loss should not go below the theoretical minimum value. However, it is surprising to see that the training loss becomes significantly lower than the dotted line at about $25000$ iterations so consequently the solution obtained afterwards is completely different from the exact one, although is not shown here. To find out what went wrong, we carefully re-evaluated all terms in the loss function (\ref{Eq:loss}) using $10^6$ testing points throughout training to ensure accurate evaluations of its continuous counterparts in the energy. The results are shown as dashed (red) lines in Fig.~\ref{Fig:fixed_point_loss}. We infer these values as the actual corresponding energy values in  training process, as they are more accurately representing the true energy values. It can also be seen in panel (a) that the actual energy during training is indeed always larger than the theoretical minimum value which is consistent with our previous analysis. Besides, at about $25000$ iterations, the actual energy is increasing, unlike the training loss that continues to decrease. We also observe that the contributions of loss from the domain boundary and interface are accurately predicted while the one from the domain (shown in panel (b)) is far from the correct value at later stage of training. Consequently the training loss is totally different from its actual energy during training (see Fig.~\ref{Fig:fixed_point_loss}(a) at later stage of iterations) and predicts the values that are much less than the correct ones.

\begin{figure}[h]
\centering
\includegraphics[scale=0.7]{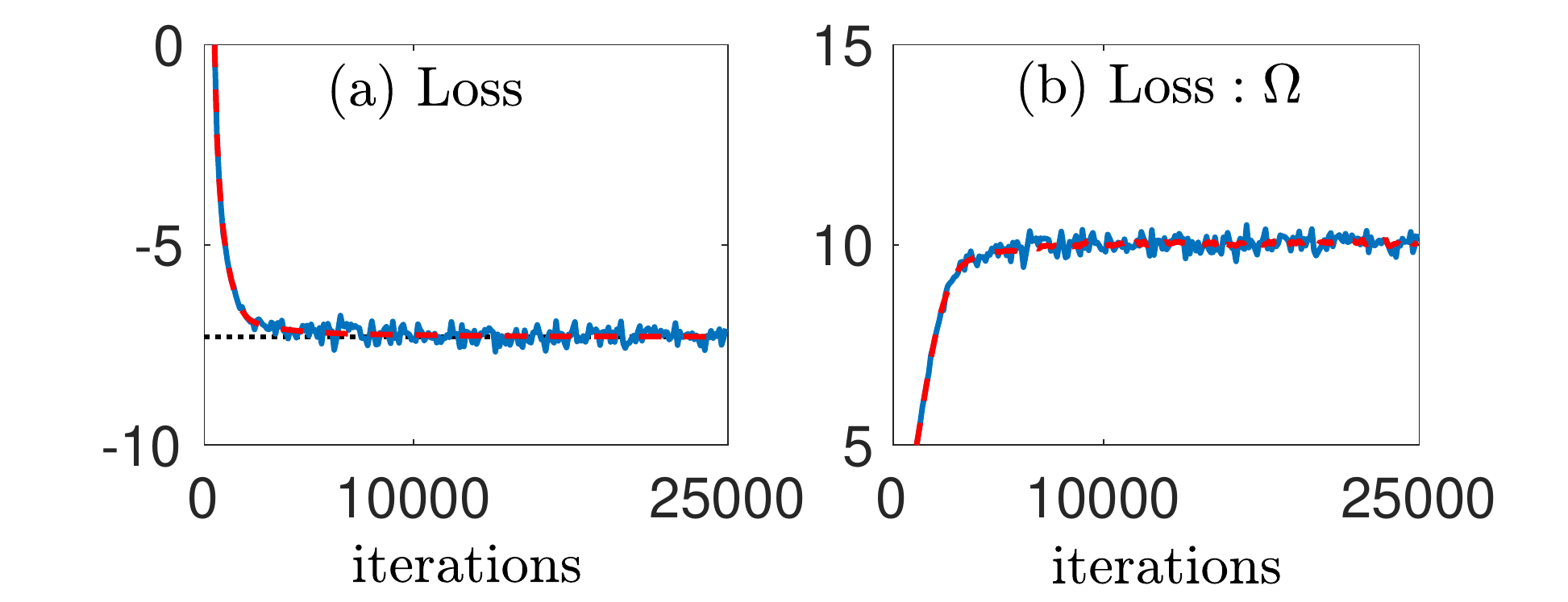}
\caption{Evaluation of the loss functions. The solid (blue) line shows the loss obtained during training process while the dashed (red) line shows a re-evaluation of the loss using $10^6$ testing points. The dotted (black) line in (a) indicates the theoretical energy value corresponding to the exact solution.}
\label{Fig:example1_2_loss}
\end{figure}

Such a scenario is known as over-fitting: the model learns too much detail of the training data and it deteriorates the performance of the model on new data. We therefore conclude that using fixed training points causes the problem of over-fitting. We should also point out that it happens even when the number of training points is much greater than the number of parameters (here in this experiment $M+M_\Gamma+M_b=1920$ vs. $N_p=151$). To overcome such a difficulty, in the following we use a strategy similar to the mini-batch gradient descent method; that is, we re-select the training points randomly in each iteration during the optimization process. As the training set is not fixed throughout the entire process, the loss evaluation in each iteration is made with a fresh set of points and therefore does not inherently overfit. The above confirmation is shown in Fig.~\ref{Fig:example1_2_loss}(a) where we solve Example 1 again, but with $1600$ randomly selected training points in each iteration. This time, the training and testing loss, shown as solid (blue) and dashed (red) lines respectively, align nicely with the theoretical energy minimum, which shows clearly the training is not over-fitted. Meanwhile, the contribution of loss from the domain also predicts correctly, see Fig.~\ref{Fig:example1_2_loss}(b).

\subparagraph{\textbf{Comparison of Optimizers}}

\begin{figure}[h]
\centering
\includegraphics[scale=0.35]{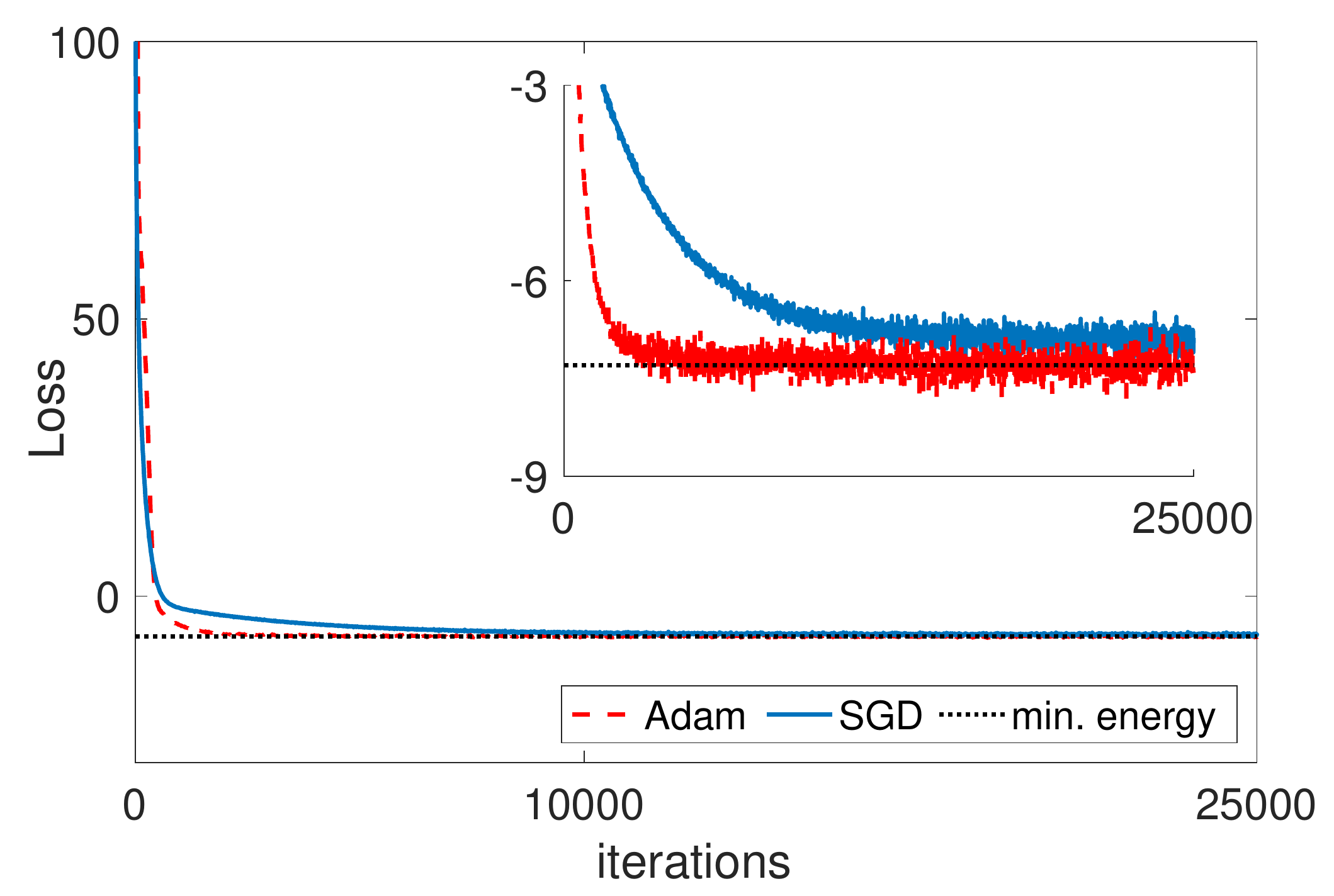}
\caption{Comparison between the performance of optimizers. The dotted (black) line indicates the theoretical minimum energy given by the exact solution. The inset shows a zoom-in to the region close to the minimum energy. }
\label{Fig:optimizer_loss}
\end{figure}

With the above re-selection strategy, we then compare the training performance with two different optimizers: Adam and stochastic gradient descent (SGD). The results are shown in Fig.~\ref{Fig:optimizer_loss} where the dashed (red) and solid (blue) lines show the training losses corresponding to Adam and SGD, respectively. The dotted (black) line indicates the theoretical minimum energy attained by the exact solution. It can be seen that the losses of the two optimizers decrease rapidly, while the Adam optimizer converges slightly faster.

We also show an enlarged view of the region near the minimum energy in the inset of the figure. The loss oscillates because of the stochastic nature of the training points selection, as expected. But on average, Adam converges faster and is closer (more accurate) to the minimum energy than SGD. Therefore, in all the following examples, we use Adam optimizer.

\subparagraph{\textbf{Comparison with the immersed boundary method}}

Here we compare the solution of the present shallow Ritz method with the numerical solution obtained by the IB method, which is known to be a first-order accurate finite difference method for equations~(\ref{eq:poisson}-\ref{eq:bdc}) on Cartesian grids~\cite{BL92} due to the employment of the discrete (or regularized) delta function. Note that in the IB method, the total number of degrees of freedom (unknowns), denoted by $N_{deg}$, is equal to the sum of the number of the Cartesian grid points $m^2$. We also denote $m_\Gamma$ as the number of Lagrangian markers on the interface.

Table~\ref{Table:exp1} shows the comparison results. The IB method uses the grid resolutions $m=m_\Gamma=80$, $160$, and $320$ corresponding to the number $N_{deg}=6400$, $25600$, and $102400$, respectively, while the present shallow Ritz method uses $N = 10$, $20$, and $30$ neurons in the hidden layer of the network, corresponding to just $N_p=51$, $101$, and $151$ parameters. One can see how significantly different those numbers of unknowns are. Using just a few number of neurons with training points $(M, M_\Gamma, M_b) = (200,80,80)$, and the penalty constant $\beta=200$, the results obtained by the present network are comparable with the ones obtained by the IB method even in the relative $L_\infty$ errors.

The accuracy of the proposed shallow Ritz method possibly depends on several factors. One is the penalty constant $\beta$. With a finite penalty constant, the minimizer of the energy functional is changed to satisfy a Robin-type boundary condition that introduces an $O(\beta^{-1})$ error. The second one is the quadrature rule that approximates the energy functional. Here, we choose the Monte Carlo type quadrature rule that has convergence rate $O(M^{-1/2})$, where $M$ is the number of sampling points. The third one is the precision of the computation.  Neural networks can be easily implemented on GPUs to take full advantage of parallel computing but at the expense of accuracy.
However, we also implement the code running in double precision by fixing $\beta$ and $M$ but does not seem to improve the overall accuracy significantly. So once the number of neurons is sufficient (here $N=20$), the magnitude of the error becomes steady. This explains why an increase in the number of neurons to $N=30$ does not help improve the accuracy.

\begin{table}[h]
\centering
\begin{tabular}{c|c||c|c|cc}
\hline
$N_{deg}$ & $\|u_{IB}-u\|_{\infty}/ \|u\|_{\infty}$  & $(N,\, N_p)$ & $\|u_\mathcal{S}-u\|_{\infty}/\|u\|_{\infty}$ & $\|u_\mathcal{S}-u\|_2/ \|u\|_2$  \\
\hline
6400 & 1.9333e$-$02  & (10,\,51) & 1.8172e$-$02 & 1.1883e$-$02 \\
25600 & 9.7151e$-$03  & (20,\,101) & 9.5521e$-$03 & 6.7409e$-$03 \\
102400 & 4.8355e$-$03 & (30,\,151) & 7.5025e$-$03 & 6.8292e$-$03  \\
\hline
\end{tabular}
\caption{Comparison between immersed boundary method and shallow Ritz method. $u_{IB}$: Solution obtained by the IB method. $u_\mathcal{S}$: Solution obtained by the present model. $u$: Exact solution. $(M, M_\Gamma, M_b) = (200,80,80)$, $\beta=200$.}
\label{Table:exp1}
\end{table}

\paragraph{\textbf{Example 2}}

In the second example, we aim to demonstrate the effectiveness of the present level set function augmentation described in Subsection~\ref{smooth} by solving the problems with or without the level set function input. We choose a setup similar to the previous example, a square domain $\Omega=[-1, 1]\times [-1, 1]$ with a circular interface of radius $0.5$ centered at the origin. We choose $\alpha=1$ and the exact solution is given as
\begin{align}
    u(x,y) =
    \left\{
    \begin{array}{ll}
        - \ln\left(x^2+y^2\right)+\sin(x)+\sin(y) & (x,y) \in\Omega^+,\\
        -\ln(0.5^2)+\sin(x)+\sin(y) & (x,y) \in \Omega^-.\\
    \end{array}\right.
\end{align}

If one does not require an augmented variable, i.e., not taking the third input of a level set function, a shallow network can be developed with just $2$ neurons in the input layer (recall for the present shallow Ritz method, there are $3$ neurons in the input layer for two-dimensional problems). For this network, the total number of parameters becomes $(d+2)N+1$ where $N$ is the number of neurons used in the hidden layer.  The network can be trained with exactly the same loss function defined in Eq.~(\ref{Eq:loss}), and we denote the solution without augmentation as $u_\mathcal{T}$.

Table~\ref{Table:exp_2input} shows a comparison between the accuracy performance of $u_\mathcal{T}$ (without level set function input) and  $u_\mathcal{S}$ (with level set function input) where the number of training points are given as $(M, M_\Gamma, M_b) = (1600, 160, 160)$. In fact, it is not required using such many training points in the domain, but we just want to avoid insufficient training points in all the following runs.

Indeed, the augmented input carrying level set function information can effectively improve training accuracy. For $u_\mathcal{T}$, the network output fails to capture the exact solution as one can see the relative $L_\infty$ error is of the order $10^{-1}$ even with $500$ neurons in the hidden layer ($2001$ parameters). However, for the present network with level set function augmentation, using $20$ neurons in the hidden layer is already capable to approximate the solution to the order $10^{-3}$ both in relative $L_\infty$ and $L_2$ errors. This example precisely shows the effectiveness of the augmented variable for carrying additional level set information to the network.

\begin{table}[h]
\centering
\begin{tabular}{c|cc||c|cc}
\hline
$(N,\, N_p)$ & $\|u_\mathcal{T}-u\|_{\infty}/ \|u\|_{\infty}$ &$\|u_\mathcal{T}-u\|_2/ \|u\|_2$ & $(N,\, N_p)$ & $\|u_\mathcal{S}-u\|_{\infty}/\|u\|_{\infty}$ & $\|u_\mathcal{S}-u\|_2/ \|u\|_2$  \\
\hline
(10,\,41) & 4.0140e$-$01 & 3.9491e$-$01 & (10,\,51) & 1.7232e$-$02 & 8.3483e$-$03 \\
(30,\,121) & 3.6079e$-$01 & 2.4959e$-$01  & (20,\,101) & 7.1492e$-$03 & 3.7503e$-$03 \\
(500,\,2001) & 3.2370e$-$01 & 2.1671e$-$01  & (30,\,151) & 5.5734e$-$03 & 3.6137e$-$03  \\
\hline
\end{tabular}
\caption{Comparison between the networks with or without the third level set function input.  $u_\mathcal{T}$: Solution obtained by a shallow network with $2$ neurons in the input layer. $u_\mathcal{S}$: Solution obtained by the present model. $u$: Exact solution. $(M, M_\Gamma, M_b) = (1600, 160, 160)$, $\beta=200$.
}
\label{Table:exp_2input}
\end{table}

To quantify further the sources of error in the proposed method, in Table~\ref{Table:exp_2beta} we show the experiments to see the effect of the penalty constant $\beta$. It reveals that the error is linearly proportional to $1/\beta$. Such a result is also consistent with the analysis in Sec.~2.2, where we have shown the global minimizer satisfies a boundary condition that is changed by $O(1/\beta)$ to the true one.

\begin{table}[h]
\centering
\begin{tabular}{c|cc}
\hline
 $\beta$  & $\|u_{\mathcal{S}}-u\|_{\infty} / \|u\|_{\infty}$ & $\|u_{\mathcal{S}}-u\|_2 / \|u\|_2 $  \\
\hline
  1   &3.6104e$-$01 & 5.3255e$-$01 \\
 10 &4.9901e$-$02 & 6.4870e$-$02 \\
 100 &6.9001e$-$03 & 6.7623e$-$03  \\
\hline
\end{tabular}
\caption{Comparison between the networks with different penalty constant $\beta$.  $u_\mathcal{S}$: Solution obtained by the present model in Example 2. $u$: Exact solution. $(M, M_\Gamma, M_b) = (1600, 160, 160)$. The number of neurons and number of parameters $(N, \, N_p)=(30, 151)$.}
\label{Table:exp_2beta}
\end{table}

We also conduct experiments to test the effect of network depth structure on the error. Here, we consider $k$-hidden layer network, where $k=1,2,3$. To have a fair comparison, we choose the number of neurons in each hidden layer $N_k$ accordingly so that the total parameters $N_p$ to be trained are similar. Table~\ref{Table:exp_2dnn} shows the relative $L_\infty$ and $L_2$ errors for those $k$-hidden layer networks. One can immediately see that for networks with roughly the same number of parameters, the magnitudes of the error are quite similar. Also, an increase in the number of parameters does not really improve the accuracy for the present problem, as clearly shown by the result of the three-hidden layer network $(k, N_k, N_p)=(3, 30, 2011)$ where the error remains in the order of $10^{-3}$.
\begin{table}[h]
\centering
\begin{tabular}{c|cc}
\hline
 $(k, N_k, N_p)$&$\|u_{\mathcal{S}}-u\|_{\infty} / \|u\|_{\infty}$ & $\|u_{\mathcal{S}}-u\|_2 / \|u\|_2 $  \\
\hline
(1, 30, 151) & 5.5734e$-$03 & 3.6137e$-$03\\
(2, 10, 161) &9.0884e$-$03 & 5.8180e$-$03 \\
(3, 7,  148)  &7.4770e$-$03 & 6.2458e$-$03 \\
\hline
(3, 30, 2011)  & 9.2629e$-$03 & 7.7047e$-$03 \\
\hline
\end{tabular}
\caption{Comparison between the networks with different depths.  $u_\mathcal{S}$: Solution obtained by the present model in Example 2. $u$: Exact solution. $(M, M_\Gamma, M_b) = (1600, 160, 160)$, $\beta=200$.}
\label{Table:exp_2dnn}
\end{table}

\paragraph{\textbf{Example 3}}

In the third example, we highlight the mesh-free nature of neural network by considering an irregular domain $\Omega$ that is given by a polar curve $r(\theta) = 1 - 0.2\cos(5\theta)$. The interface $\Gamma$ is chosen as an ellipse that can be labeled by the zero level set of the function $\phi(x, y) = \frac{x^2}{0.7^2}+\frac{y^2}{0.5^2}-1$. We fix $\alpha=0$ and the exact solution is chosen as
\begin{align}
    u(x,y) =
    \left\{
    \begin{array}{ll}
        \ln\left(\frac{x^2}{0.7^2}+\frac{y^2}{0.5^2}\right) & (x,y) \in\Omega^+,\\
        \sin(x)\cos(y)\left[\left(\frac{x^2}{0.7^2}+\frac{y^2}{0.5^2}\right)^2-1\right]  & (x,y) \in \Omega^-.\\
    \end{array}\right.
\end{align}

We randomly sample training points in the domain with $(M, M_\Gamma, M_b) = (400, 80, 80)$. The results are shown in Table~\ref{table:2D_irregular} and the solutions are accurate to the order of $10^{-3}$ in relative $L_2$ error with just $30$ neurons in the hidden layer. We should also point out that it can be quite tedious to implement for traditional finite difference methods to solve the problem on such irregular domain. Thus, we emphasize that the irregular domain case can be properly handled with no substantial difficulty.

\begin{table}[h]
    \centering
    \begin{tabular}{c|cc}
    \hline
    $(N,\, N_p)$ & $\|u_\mathcal{S}-u\|_{\infty}/\|u\|_{\infty}$ & $\|u_\mathcal{S}-u\|_2/\|u\|_2$  \\
    \hline
    (30, 151) &  1.7095e$-$02 & 8.8148e$-$03\\
    (40, 201) &  1.2277e$-$02 & 8.5489e$-$03\\
    \hline
    \end{tabular}
    \caption{$u_\mathcal{S}$: Solution obtained by the present model in Example 3. $u$: Exact solution. $(M, M_\Gamma, M_b) = (400,80,80)$, $\beta=200$.}
    \label{table:2D_irregular}
\end{table}

\begin{figure}[h]
\centering
\includegraphics[scale=0.5]{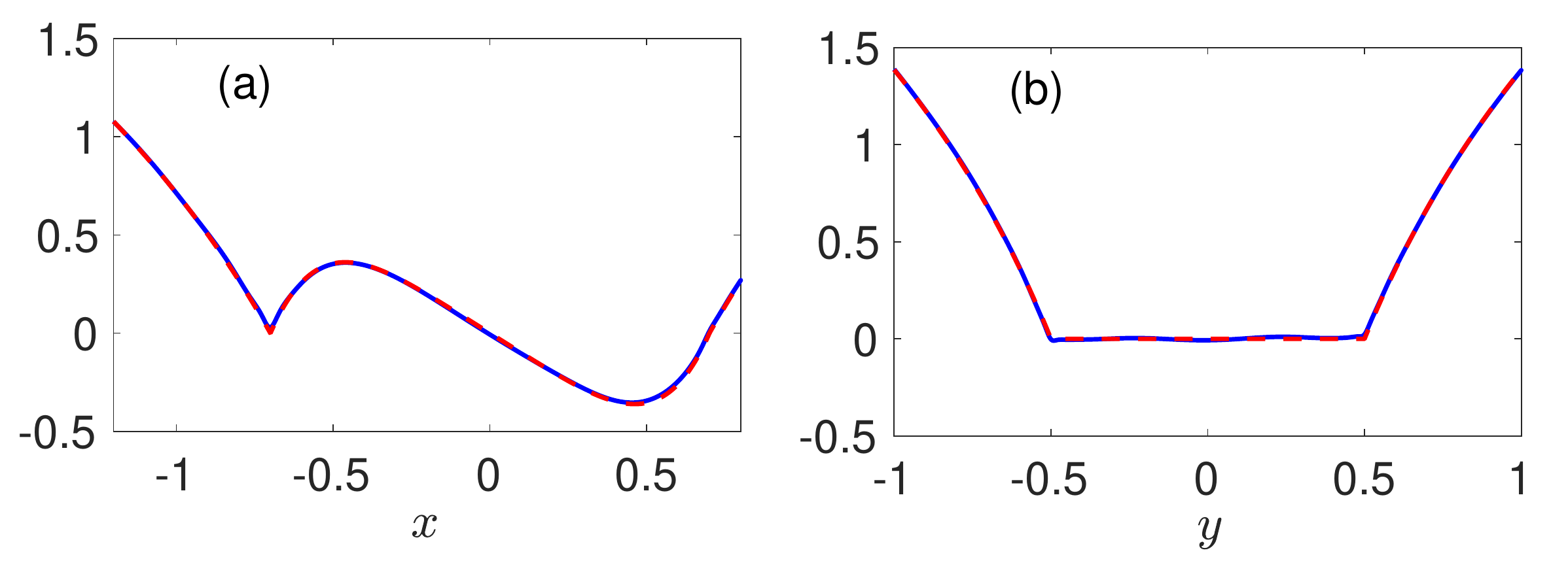}
\caption{Cross sections at (a) $y=0$ and (b) $x=0$ for the solution in Example 3 with $N=40$. The model solution $u_\mathcal{S}$ is shown as solid (blue) line and the exact solution is shown as dashed (red) line.}
\label{Fig:2D_irregular_cross}
\end{figure}

\begin{figure}[h!]
\centering
\includegraphics[scale=0.3]{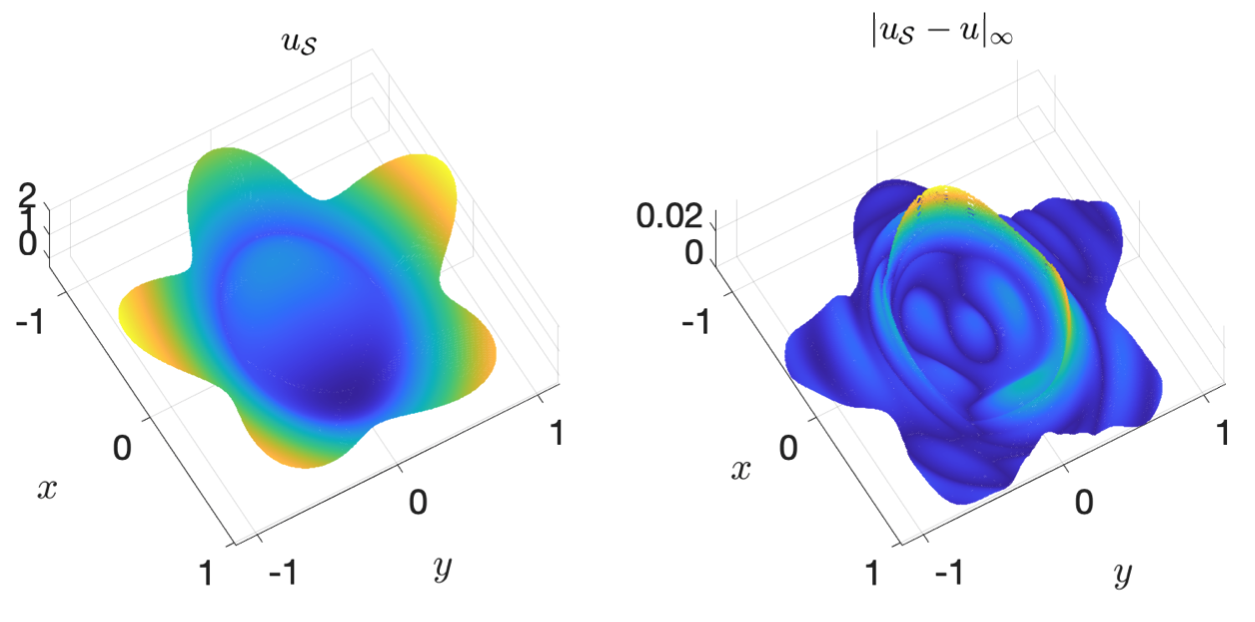}
\caption{Left: The solution profile obtained by the present shallow Ritz method with $N=40$. Right: The profile of absolute error.}
\label{Fig:2D_irregular_profile}
\end{figure}

We also show the cross-sections of the solution at $y=0$ and $x=0$ in Figure~\ref{Fig:2D_irregular_cross}(a) and (b), respectively.  One can see that, the present network indeed accurately approximates the solution even though there are cusps at the interface. Furthermore, in Fig.~\ref{Fig:2D_irregular_profile} showing the solution profile and absolute error correspondingly, we observe that the largest error occurs in the region close to the interface, but is not significantly larger than the values in other regions.

\paragraph{\textbf{Example 4}}

In this example, we show the ability of the present method to solve three-dimensional problems. We choose the domain as a cube $\Omega = [-1,1]^3$ with a spherical interface that is labeled by a zero-level set of the function $\phi(x,y,z) = \frac{x^2}{0.4^2} + \frac{y^2}{0.4^2}+\frac{z^2}{0.4^2}-1$. The exact solution is chosen as
\begin{align}
    u(x,y,z) =
    \left\{
    \begin{array}{ll}
        x(-1 + \exp(0.4^2-(x^2 + y^2+z^2))) & (x,y,z) \in\Omega^+,\\
        -1 + \cos(0.4^2-(x^2 + y^2+z^2))  & (x,y,z) \in \Omega^-.\\
    \end{array}\right.
\end{align}
In this test, we choose $\alpha=1$ and $(M,M_\Gamma,M_b) = (216,216,216)$. The results are shown in Table~\ref{table:3D}. Again, the present network shows a good performance. The relative $L_2$ error is less than $1\%$ with just $30$ neurons in the hidden layer.

\begin{table}[!]
    \centering
    \begin{tabular}{c|cc}
    \hline
    $(N,\, N_p)$ & $\|u_\mathcal{S}-u\|_{\infty}/\|u\|_{\infty}$ & $\|u_\mathcal{S}-u\|_2/\|u\|_2$  \\
    \hline
    (20, 121) &  1.9960e$-$02 & 1.4343e$-$02\\
    (30, 181) &  1.6274e$-$02 & 9.9769e$-$03\\
    (40, 241) &  1.2727e$-$02 & 8.7665e$-$03\\
    \hline
    \end{tabular}
    \caption{$u_\mathcal{S}$: Solution obtained by the present model in Example 4 (three-dimensional case). $u$: Exact solution. $(M, M_\Gamma, M_b) = (216, 216, 216)$, $\beta=100$.}
    \label{table:3D}
\end{table}

\paragraph{\textbf{Example 5}}

As the last example, we demonstrate the capability of the present method to solve high-dimensional problem by taking the dimension size $d=6$. For the problem setup, we consider a $6$-sphere of radius $0.6$ as the domain $\Omega$ enclosing another smaller $6$-sphere of radius $0.5$ as the interface $\Gamma$. The interface $\Gamma$ can be labeled by the zero level set of the function $\phi(\bx) = \left(\frac{\|\mathbf{x}\|_2}{0.5}\right)^2- 1$. We fix $\alpha=0$ and the exact solution is chosen as
\begin{align}
    u(\mathbf{x}) =
    \left\{
    \begin{array}{ll}
        \exp(0.5^2-\|\mathbf{x}\|_2^2) + \sum^5_{i=1}\sin(x_i)& \mathbf{x} \in\Omega^+,\\
        1+2\sin(0.5^2 - \|\mathbf{x}\|_2^2) + \sum^5_{i=1}\sin(x_i) & \mathbf{x} \in \Omega^-,\\
    \end{array}\right.
\end{align}
where $\mathbf{x} = (x_1, x_2,x_3,x_4,x_5,x_6)$.

Here, we choose the number of training points based on the following strategy. Given a radius $R$, the volume and surface area of the $6$-sphere are $(1/6)\pi^3R^6$ and $\pi^3R^5$, respectively, so the ratio between these two numbers is $R^6 : 6R^5$. We then choose the number of training points in the sphere and on the surface based on this ratio; that is, if we have $500$ points in the domain that corresponds to the effective radius $R= (500)^{1/6}$, then we select $6R^{5} \approx 1065$ points on the boundary and interface.

\begin{figure}[h!]
\centering
\includegraphics[scale=0.5]{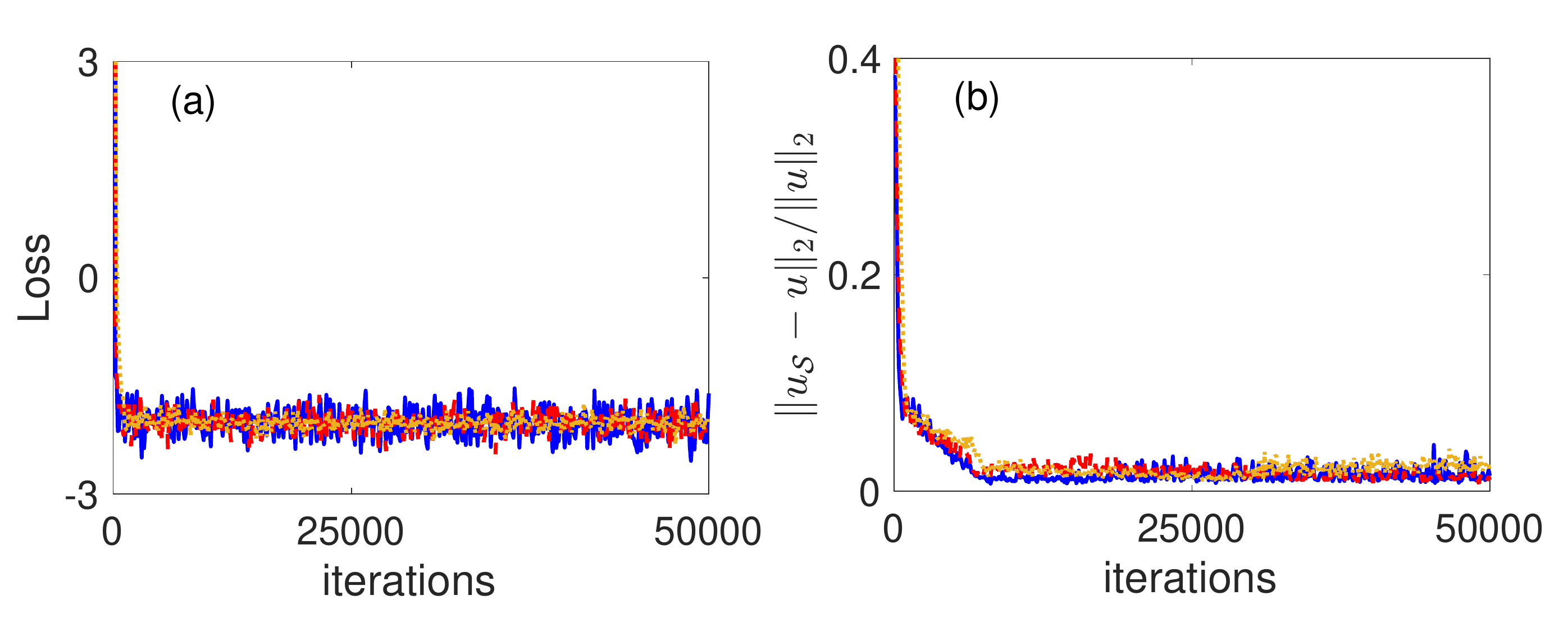}
\caption{ The evolution plots of (a) training loss, and (b) relative $L_2$ error. $(M, M_\Gamma, M_b)=(100, 278, 278)$: solid (blue) line; $(200, 496, 496)$: dashed (red) line; $(500, 1065, 1065)$: dotted (yellow) line.}
\label{Fig:6D}
\end{figure}

The results are shown in Table~\ref{table:6D}. With just $10$ neurons in the hidden layer, the relative $L_2$ error is in the order of $10^{-3}$. In addition, we compare the performance between different number of training points, shown in Table~\ref{table:6D_2}. It is surprising to see that accuracy of the solution is already good enough by choosing $(M, M_\Gamma, M_b) = (100, 278, 278)$ in six dimensions, while doubling or even quintupling the training points does not improve the accuracy. More precisely, we also show  the evolution of the training loss and the relative $L_2$ error during the training process in Fig.~\ref{Fig:6D}. One can see that the results are almost indistinguishable from three different training point choices. Here, even in $6$ dimensions, we only need as few as $656$ training points to achieve the desired accuracy. A possible explanation is that the error of the Monte Carlo integration is independent of the problem dimension. Although we have conducted other possible sources of error such as the choice of penalty number $\beta$ and the network depth in Example 2, a detailed error analysis is still required for further investigation.

\begin{table}[h]
    \centering
    \begin{tabular}{c|cc}
    \hline
    $(N,\, N_p)$ & $\|u_\mathcal{S}-u\|_{\infty}/\|u\|_{\infty}$ & $\|u_\mathcal{S}-u\|_2/\|u\|_2$  \\
    \hline
    (10, 91) & 2.8309e$-$02    & 6.4073e$-$03 \\
    (20, 181) & 2.6612e$-$02    & 7.0114e$-$03    \\
    (30, 271) & 2.0292e$-$02& 6.9735e$-$03  \\
    \hline
    \end{tabular}
    \caption{$u_\mathcal{S}$: Solution obtained by the present model in Example 5 (six-dimensional case). $u$: Exact solution. $(M, M_\Gamma, M_b) = (500, 1065, 1065)$, $\beta=100$.}
    \label{table:6D}
\end{table}

\begin{table}[h]
    \centering
    \begin{tabular}{c|c|cc}
    \hline
    $(M, M_\Gamma, M_b)$ & $(N,\, N_p)$ & $\|u_\mathcal{S}-u\|_{\infty}/\|u\|_{\infty}$ & $\|u_\mathcal{S}-u\|_2/\|u\|_2$  \\
    \hline
 $(100, 278, 278)$ &  (10, 91) & 2.4877e$-$02    & 7.4379e$-$03    \\
 $(200, 496, 496)$  & (10, 91) & 2.5073e$-$02  & 7.2977e$-$03    \\
  $(500, 1065, 1065)$  & (10, 91) & 2.8309e$-$02    & 6.4073e$-$03    \\
    \hline
    \end{tabular}
    \caption{Comparison between the number of training points in Example 5. $u_\mathcal{S}$: Solution obtained by the present model. $u$: Exact solution.}
    \label{table:6D_2}
\end{table}

\section{Conclusion}

In this paper, a novel shallow Ritz method is developed to solve elliptic problems with delta function singular sources on the interface. By introducing an energy functional, we reformulate the governing equation as a variational problem. The crucial observation is that the contribution of the singular source in the equation becomes a regular surface integral term in the energy. Therefore, we do not need to introduce the discrete delta function used in traditional regularization methods, such as the immersed boundary method. To enforce the boundary condition, we add a penalty term to the energy and find that this treatment changes the global minimizer to one that satisfies Robin-type boundary condition.

We propose a shallow neural network to approximate the global minimizer of the energy functional, and train the network by minimizing a loss function that presents a discrete version of the energy. In addition, we include the level set function as an additional feature input to the network and find that it significantly improves the training accuracy. We perform a series of numerical tests to show the accuracy and efficiency of the present network and its capability to handle problems in irregular domains or high dimensions. As shown in numerical experiments in this paper, most testing problems can be solved with acceptable accuracy by the present network with moderate number of neurons (no larger than $40$).
Although the present network is similar in spirit to the deep Ritz method~\cite{EY18}, here we consider a completely shallow one (only one hidden layer) so it significantly reduces the computational complexity and learning workload without sacrificing the accuracy.

We have to emphasize that it is not our intention to compete with the traditional methods such as immersed boundary (IB) method, immersed interface method (IIM), immersed finite element (IFE) method or other grid-based methods listed in the reference. Instead, we just provide an alternative using neural network method and follow the recent pioneering deep Ritz method by E and Yu in \cite{EY18} to solve the elliptic interface problems with singular delta function sources. Since the underlying solution of the problem is usually not smooth (the function is continuous but not for its derivatives) across the interface, the traditional accurate methods (for instance IIM, IFE method, see the book of Li and Ito [24]) need special treatments near the interface. However, due to the mesh-free advantage, the present neural network uses the level set function of the interface as a feature input and learns the solution directly.  Once the network is trained successfully, the solution at any point in the domain can be obtained by a feedforward one hidden layer neural network approximation. So from this point of view, it is efficient. In fact, the training time of the present network takes only a few minutes on a desktop, even for a 6-dimensional problem listed in Example 5. Meanwhile, neural networks can be easily implemented on GPUs to take full advantage of parallel computation and are easy to deal with problems defined in irregular domains as well. Another advantage is that the present approach can handle high dimensional problems with exactly same framework and the number of hyperparameters (unknowns in neural network method) to be trained scales only linearly with the dimension.

As mentioned before, the accuracy of the proposed shallow Ritz method possibly depends on three factors, namely the penalty constant $\beta$, the quadrature rule for estimating the energy functional, and the precision of the computation. To improve accuracy, we must consider three factors simultaneously. We believe that the present one hidden layer network is not the reason caused the accuracy no better than $10^{-3}$ in relative $L_\infty$ or $L_2$ error throughout all numerical experiments here. The one hidden layer network with sufficiently smooth activation function has been proved to be capable to represent an arbitrary function and its derivatives to arbitrary accuracy. From the numerical results shown in Section 4, the accuracy obtained by the present method is comparable with the one existing in related literatures despite the fact that the present network has simplest structure (one hidden layer).   For example, the pioneering work of E and Yu in \cite{EY18} where the authors used the deep learning network ResNet (say a stack of 4 blocks with 2 fully-connected layers in each block), the relative $L_2$ errors for all considered problems are also up to $10^{-3}$. In the recent work of Guo and Yang in \cite{GY21}, a similar ResNet was applied to solve elliptic interface problems where the jump conditions are put into their unfitted Nitsche's energy functional. Again, most of their numerical results in  relative $L_2$ errors are no better than $10^{-3}$.
Certainly, it will be interesting to investigate further to improve the accuracy of the Ritz-type neural network method for solving PDEs which we shall leave as our future work. Nevertheless, from the authors' point of view, the Ritz-type neural network is still valuable and useful for computational physics applications due to following two reasons. First, it is completely mesh-free for high-dimensional problems when traditional numerical methods are hard to tackle. Second, in many interesting physical problems, finding a solution to the problem is often equivalent to finding the minimum of its energy law, which falls into the Ritz-type neural network methodology naturally.

The present goal is to solve some elliptic interface problems using a Ritz-type neural network. Indeed, for given $\Omega$, $\Gamma$, $f(\bx)$, $c(\bs)$ and the boundary condition, we need to train a new network that satisfies the equation and boundary condition simultaneously. These given information are all from the original PDE formulation so there is no difficulty in generating the training data. (That is, generating training data from a well-posed PDE can be done without knowing any numerical techniques for solving that PDE.) However, when using some traditional numerical techniques such as IB method for the present problem, if the domain $\Omega$ and mesh size are fixed, the matrix of resultant linear system remains the same even if the interface $\Gamma$ and the right-hand side functions $f(\bx)$, $c(\bs)$ are changed. So from this point of view, the IB method is more generalizable. Therefore, to improve the generalizability of the present network,
we might consider to design a neural network as an operator representation that has the solution as the output. But this is beyond the scope of the paper and we leave it as our future work.

In the present work, we only consider stationary elliptic problems with constant coefficients. There should be no difficulty in considering problems with contrast coefficients or even variable coefficients. The same framework can be applied straightforwardly by writing the corresponding energy functionals. As a forthcoming extension, we shall consider time-dependent problems, and particularly the moving interface problems.

\section*{Acknowledgement}

T.-S. Lin and W.-F. Hu acknowledge supports by Ministry of Science and Technology (MOST), Taiwan, under research grant 109-2115-M-009-006-MY2 and 109-2115-M-008-014-MY2, respectively.
T.-S. Lin and W.-F. Hu also acknowledge support by NCTS of Taiwan.
M.-C. Lai acknowledges the support by MOST, Taiwan, under research grants 108-2119-M-009-012-MY2 and 110-2115-M-A49-011-MY3.



\end{document}